\RequirePackage{ifpdf}
\ifpdf 
\documentclass[pdftex]{sigma}
\else
\documentclass{sigma}
\fi

\begin{document}

\renewcommand{\PaperNumber}{006}

\FirstPageHeading

\ShortArticleName{Peterson's Deformations of Higher Dimensional Quadrics}

\ArticleName{Peterson's Deformations\\ of Higher Dimensional Quadrics}

\Author{Ion I. DINC\u{A}}

\AuthorNameForHeading{I. Dinc\u{a}}

\Address{Faculty of Mathematics and Informatics,
University of Bucharest, \\ 14 Academiei Str., 010014, Bucharest,
Romania}
\Email{\href{mailto:dinca@gta.math.unibuc.ro}{dinca@gta.math.unibuc.ro}}

\ArticleDates{Received July 13, 2009, in f\/inal form January 16, 2010;  Published online January 20, 2010}

\Abstract{We provide the f\/irst explicit examples of deformations
of higher dimensional quadrics: a straightforward generalization
of Peterson's explicit $1$-dimensional family of deformations in
$\mathbb{C}^3$ of $2$-dimensional general quadrics with common
conjugate system given by the spherical coordinates on the complex
sphere $\mathbb{S}^2\subset\mathbb{C}^3$ to an explicit
$(n-1)$-dimensional family of deformations in $\mathbb{C}^{2n-1}$
of $n$-dimensional general quadrics with common conjugate system
given by the spherical coordinates on the complex sphere
$\mathbb{S}^n\subset\mathbb{C}^{n+1}$ and non-degenerate joined
second fundamental forms. It is then proven that this family is
maximal.}


\Classification{53A07; 53B25; 35Q58}

\section{Introduction}

The Russian mathematician Peterson was a student of Minding's, who
in turn was interested in deformations (through bending) of
surfaces\footnote{See Peterson's biography at \url{http://www-history.mcs.st-and.ac.uk/Biographies/Peterson.html}.}, but unfortunately most of his works
(including his independent discovery of the Codazzi--Mainardi
equations and of the Gau\ss--Bonnet theorem) were made known to
Western Europe mainly after they were translated in 1905 from
Russian to French (as is the case with his deformations of
quadrics~\cite{P}, originally published in 1883 in Russian).
Peterson's work on deformations of general quadrics preceded that
of Bianchi, Calapso, Darboux, Guichard and \c{T}i\c{t}eica's from
the years 1899--1906 by two decades; in particu\-lar Peterson's
$1$-dimensional family of deformations of surfaces admitting a~common {\it conjugate system}~$(u,v)$ (that is the second
fundamental form is missing mixed terms $du\odot dv$) are {\it
associates} (a~notion naturally appearing in the inf\/initesimal
deformation problem) to Bianchi's $1$-dimensional family of
surfaces satisfying $(\log K)_{uv}=0$ in the common asymptotic
coordinates~$(u,v)$, $K$~being the Gau\ss\ curvature (see Bianchi~\cite[Vol.~2, \S\S~294, 295]{B1}).

The work of these illustrious geometers on deformations in
$\mathbb{C}^3$ of quadrics in $\mathbb{C}^3$ (there is no other
class of surfaces for which an interesting theory of deformation
has been built) is one of the crowning achievements of the golden
age of classical geometry of surfaces and at the same time it
opened new areas of research (such as af\/f\/ine and projective
dif\/ferential geometry) continued later by other illustrious
geometers (Blaschke, Cartan, etc.).

Peterson's $1$-dimensional family of deformations of
$2$-dimensional quadrics is obtained by imposing an ansatz
naturally appearing from a geometric point of view, namely the
constraint that the common conjugate system of curves is given by
intersection with planes through an axis and tangent cones
centered on that axis; thus this result of Koenigs (see Darboux
\cite[\S\S~91--101]{D}) was again (at least when the cones are
tangent along plane curves) previously known to Peterson. Note
also that Calapso in \cite{Ca} has put Bianchi's B\"{a}cklund
transformation of deformations in~$\mathbb{C}^3$ of general
$2$-dimensional quadrics with center in terms of common conjugate
systems (the condition that the conjugate system on a
$2$-dimensional quadric is a conjugate system on one of its
deformations in~$\mathbb{C}^3$ was known to Calapso for a decade,
but the B\"{a}cklund transformation for general quadrics eluded
Calapso since the common conjugate system was best suited for this
transformation only at the analytic level).

Although this is the original approach Peterson used to f\/ind his
deformations of quadrics, other features of his approach will make
it amenable to higher dimensional generalizations, namely the
warping of linear element (the warping of the linear element of a
plane curve to get the linear element of a surface of revolution
$(d(f\cos(u^1)))^2+(d(f\sin(u^1)))^2=(df)^2+f^2(du^1)^2$ for
$f=f(u^2)$ is such an example) and separation of variables;
post-priori the common conjugate system property may be given a
geometric explanation analogous to that in dimension $3$.

In 1919--1920 Cartan has shown in~\cite{C} (using mostly projective
arguments and his exterior dif\/ferential systems in involution and
exteriorly orthogonal forms tools) that space forms of dimension
$n$ admit rich families of deformations (depending on $n(n-1)$
functions of one variable) in surrounding $(2n-1)$-dimensional
space forms, that such deformations admit lines of curvature
(given by a canonical form of exteriorly orthogonal forms; thus
they have f\/lat normal bundle; since the lines of curvature on
$n$-dimensional space forms (when they are considered by
def\/inition as quadrics in surrounding $(n+1)$-dimensional space
forms) are undetermined, the lines of curvature on the deformation
and their corresponding curves on the quadric provide the common
conjugate system) and that the co-dimension $(n-1)$ cannot be
lowered without obtaining rigidity as the deformation being the
def\/ining quadric.

In 1983 Berger, Bryant and Grif\/f\/iths~\cite{B} proved (including by
use of tools from algebraic geometry) in particular that Cartan's
essentially projective arguments (including the exterior part of
his exteriorly orthogonal forms tool) can be used to generalize
his results to $n$-dimensional general quadrics with positive
def\/inite linear element (thus they can appear as quadrics in
$\mathbb{R}^{n+1}$ or as space-like quadrics in
$\mathbb{R}^n\times(i\mathbb{R})$) admitting rich families of
deformations (depending on $n(n-1)$ functions of one variable) in
surrounding Euclidean space $\mathbb{R}^{2n-1}$, that the
co-dimension $(n-1)$ cannot be lowered without obtaining rigidity
as the deformation being the def\/ining quadric and that quadrics
are the only Riemannian $n$-dimensional manifolds that admit a
family of deformations in $\mathbb{R}^{2n-1}$ as rich as possible
for which the exteriorly orthogonal forms tool (naturally
appearing from the Gau\ss\ equations) can be applied.

Although Berger, Bryant and Grif\/f\/iths~\cite{B} do not explicitly
state the common conjugate system property (which together with
the non-degenerate joined second fundamental forms assumption
provides a tool similar to the canonical form of exteriorly
orthogonal forms), this will turn out to be the correct tool of
dif\/ferential geometry needed to attack the deformation problem for
higher dimensional quadrics; also at least for diagonal quadrics
without center Peterson's deformations of higher dimensional
quadrics will turn out to be amenable to explicit computations of
their B\"{a}cklund transformation\footnote{See Dinc\u{a} I.I., Bianchi's
B\"{a}cklund transformation for higher dimensional quadrics,
\href{http://arxiv.org/abs/0808.2007}{arXiv:0808.2007}.}.

All computations are local and assumed to be valid on their open
domain of validity without further details; all functions have the
assumed order of dif\/ferentiability and are assumed to be
invertible, non-zero, etc when required (for all practical
purposes we can assume all functions to be analytic).

Here we have the two main theorems concerning the
$(n-1)$-dimensional family of deformations of higher dimensional
general quadrics and respectively its maximality:

\begin{theorem}\label{theorem1}
The quadric
\[
\sum_{j=0}^n\frac{(x_j)^2}{a_j}=1,\qquad
a_j\in\mathbb{C}^*
\] distinct parameterized with the conjugate
system $(u^1,\dots ,u^n)\subset\mathbb{C}^n$ given by the spherical
coordinates on the unit sphere
$\mathbb{S}^n\subset\mathbb{C}^{n+1}$:
\begin{gather*}
\mathcal{X}=\sqrt{a_0}\mathbf{C}_0e_0+
\sum_{k=1}^n\sqrt{a_k}\mathbf{C}_k\sin\big(u^k\big)e_k,\qquad
\mathbf{C}_k:=\prod_{j=k+1}^n\cos(u^j)
\end{gather*}
and the sub-manifold
\begin{gather*}
\mathcal{X}_{\mathbf{z}}=\sum_{k=1}^{n-1}\mathbf{C}_kf_k\big(\mathbf{z},u^k\big)
\big(\cos\big(g_k\big(\mathbf{z},u^k\big)\big)e_{2k-2}
+\sin\big(g_k\big(\mathbf{z},u^k\big)\big)e_{2k-1}\big)+h\big(\mathbf{z},u^n\big)e_{2n-2}
\end{gather*}
 of
$\mathbb{C}^{2n-1}$ depending on the parameters
$\mathbf{z}=(z_1,z_2,\dots ,z_{n-1})\in\mathbb{C}^{n-1}$, $z_0:=1$
and with
\begin{gather}
f_k\big(z_{k-1},z_k,u^k\big):=\sqrt{(z_{k-1}-z_k)a_0+(a_k-z_{k-1}a_0)\sin^2(u^k)},\qquad k=1,\dots ,n-1,\nonumber\\
g_k\big(z_{k-1},z_k,u^k\big):=\int_0^{u^k}\frac{\sqrt{(z_{k-1}-z_k)a_0a_k+(a_k-z_{k-1}a_0)z_ka_0\sin^2(t)}}
{(z_{k-1}-z_k)a_0+(a_k-z_{k-1}a_0)\sin^2(t)}dt,\nonumber\\
h\big(z_{n-1},u^n\big):=\int_0^{u^n}\sqrt{a_n-(a_n-z_{n-1}a_0)\sin^2(t)}dt\label{eq:fgn}
\end{gather}
have the same linear element
$|d\mathcal{X}|^2{=}|d\mathcal{X}_{\mathbf{z}}|^2$. For
$z_1={\cdots}=z_{n-1}=0$ we get \mbox{$g_2={\cdots}=g_{n-1}=0$},
$\mathcal{X}=\mathcal{X}_{\mathbf{0}}$ with
$\mathbb{C}^{n+1}\hookrightarrow\mathbb{C}^{2n-1}$ as
$(x_0,x_1,\dots ,x_n)\mapsto(x_0,x_1,x_2,0,x_3,0,\dots ,x_{n-1},0,x_n)$.
For $z_1=\dots =z_{n-1}=1$ we get
$\mathcal{X}_{\mathbf{1}}=(\mathbf{x}_0,\dots ,\mathbf{x}_{2n-2})$
given by Peterson's formulae
\begin{gather}
\sqrt{(\mathbf{x}_{2k-2})^2+(\mathbf{x}_{2k-1})^2}=\sqrt{a_k-a_0}\mathbf{C}_k\sin\big(u^k\big),\nonumber\\
 \tan^{-1}\left(\frac{\mathbf{x}_{2k-1}}{\mathbf{x}_{2k-2}}\right)
=\frac{\sqrt{a_0}}{\sqrt{a_k-a_0}}\tanh^{-1}\big(\cos\big(u^k\big)\big), \qquad k=1,\dots,n-1,
\nonumber\\
\mathbf{x}_{2n-2}=\int_0^{u^n}\sqrt{a_n-(a_n-a_0)\sin^2(t)}dt.\label{eq:pfn}
\end{gather}
Moreover $(u^1,\dots ,u^n)$ form a conjugate system on
$\mathcal{X}_{\mathbf{z}}$ with non-degenerate joined second
fundamental forms $($that is $[d^2\mathcal{X}^TN\ \
d^2\mathcal{X}_{\mathbf{z}}^TN_{\mathbf{z}}]$ is a symmetric
quadratic $\mathbb{C}^n$-valued form which contains only
$(du^j)^2$ terms for $N$ normal field of $\mathcal{X}$ and
$N_{\mathbf{z}}=[N_1\   \dots   \ N_{n-1}]$ normal frame of~$\mathcal{X}_{\mathbf{z}}$ and the dimension $n$ cannot be lowered
for $\mathbf{z}$ in an open dense set$)$.
\end{theorem}

\begin{theorem}\label{theorem2}
For $x\subset\mathbb{C}^{2n-1}$ deformation of the quadric
$x_0\subset\mathbb{C}^{n+1}$ $($that is $|dx|^2=|dx_0|^2)$ with
$n\ge 3$, $(u^1,\dots ,u^n)$ common conjugate system and
non-degenerate joined second fundamental forms,
$N_0^Td^2x_0=:\sum\limits_{j=1}^nh_j^0(du^j)^2$ second fundamental form
of $x_0$ we have $\Gamma_{jk}^l=0$ for $j$, $k$, $l$ distinct and such
deformations are in bijective correspondence with solutions
$\{\mathbf{a}_j\}_{j=1,\dots ,n}\subset\mathbb{C}^*$ of the
differential system
$\partial_{u^k}\log\mathbf{a}_j=\Gamma_{jk}^j$, $j\neq k$,
$\sum\limits_{j=1}^n\frac{(h_j^0)^2}{\mathbf{a}_j^2}+1=0.$ In particular
this implies that for $(u^1,\dots ,u^n)$ being the conjugate system
given by spherical coordinates on
$\mathbb{S}^n\subset\mathbb{C}^{n+1}$ the above explicit
$(n-1)$-dimensional family of deformations
$\mathcal{X}_{\mathbf{z}}$ is maximal.
\end{theorem}

The remaining part of this paper is organized as follows: in Section~\ref{sec:pet} we shall recall Peterson's deformations of quadrics;
the proof of Theorem~\ref{theorem1}  appears in Sections~\ref{sec:pet1},~\ref{sec:pet2} and the proof of
Theorem~\ref{theorem2} appears in Sections~\ref{sec:pet3},~\ref{sec:pet4}.

\section{Peterson's deformations of quadrics}\label{sec:pet}

Although Peterson \cite{P} discusses all types of quadrics in the
complexif\/ied Euclidean space
\[
\big(\mathbb{C}^3,\langle \cdot,\cdot\rangle\big),\qquad \langle x,y\rangle:=x^Ty,\qquad |x|^2:=x^Tx\quad \mathrm{for}\ \
x,y\in\mathbb{C}^3
\]
and their totally real cases, we shall only discuss quadrics of
the type $\sum\limits_{j=0}^2\frac{(x_j)^2}{a_j}=1$, $a_j\in\mathbb{C}^*$
distinct, since the remaining cases of quadrics should follow by
similar computations. Their totally real cases (that is
$(x_j)^2$, $a_j\in\mathbb{R}$) are discussed in detail in Peterson~\cite{P}, so we shall not insist on this aspect.

\begin{remark}
It is less known since the classical times that there are many
types of quadrics from a~complex metric point of view, each coming
with its own totally real cases (real valued (in)def\/inite linear
element); among these quadrics for example the quadric
$(x_0-ix_1)x_2-(x_0+ix_1)=0$ is rigidly {\it applicable}
(isometric) to all quadrics of its confocal family and to all its
homothetic quadrics. It is Peterson who f\/irst introduced the idea
of {\it ideal applicability} (for example a real surface may be
applicable to a totally real space-like surface
$\subset\mathbb{R}^2\times(i\mathbb{R})$ of a~complexif\/ied real
ellipsoid, so it is ideally applicable on the real ellipsoid).
\end{remark}

With $\{e_j\}_{j=0,1,2}$, $e_j^Te_k=\delta_{jk}$ the standard basis
of $\mathbb{C}^3$  and the functions $f=f(z,u^1)$, $g=g(z,u^1)$,
$h=h(z,u^2)$ depending on the parameter(s) $z=(z_1,z_2,\dots )$ to be
determined later we have the surfaces
\begin{gather}
\mathcal{X}_z:=\cos\big(u^2\big)f\big(z,u^1\big)\big(\cos\big(g\big(z,u^1\big)\big)e_0+\sin
\big(g\big(z,u^1\big)\big)e_1\big)+h\big(z,u^2\big)e_2.
\label{eq:surf}
\end{gather}
Note that the f\/ields
$\partial_{u^1}\mathcal{X}_z|_{u^1=\mathrm{const}}$,
$\partial_{u^2}\mathcal{X}_z|_{u^2=\mathrm{const}}$ generate developables
(cylinders with gene\-ra\-tors perpendicular on the third axis and
cones centered on the third axis), so $(u^1,u^2)$ is a~conjugate
system on $\mathcal{X}_z$ for every $z$; in fact all surfaces have
conjugate systems arising this way and can be parameterized as
\[
x=f\big(u^1,u^2\big)\big(\cos\big(u^1\big)e_0+\sin\big(u^1\big)e_1\big)+g\big(u^1,u^2\big)e_2, \qquad
\partial_{u^1}\left(\partial_{u^2}\left(\frac{g}{f}\right){\Big/}\partial_{u^2}\left(\frac{1}{f}\right)\right)=0.
\]

The quadric $\sum\limits_{j=0}^2\frac{(x_j)^2}{a_j}=1$ is parameterized
by the spherical coordinates
\[
\mathcal{X}=\sqrt{a_0}\cos\big(u^2\big)\cos\big(u^1\big)e_0
+\sqrt{a_1}\cos\big(u^2\big)\sin\big(u^1\big)e_1+\sqrt{a_2}\sin\big(u^2\big)e_2.
\]
We have
\begin{gather*}
|d\mathcal{X}_z|^2=\cos^2\big(u^2\big)\big(f'^2\big(z,u^1\big)
+f^2\big(z,u^1\big)g'^2\big(z,u^1\big)\big)\big(du^1\big)^2
\\
\phantom{|d\mathcal{X}_z|^2=}{}
+\tfrac{1}{2}d\big(\cos^2\big(u^2\big)\big)d\big(f^2\big(z,u^1\big)\big)
+\big(f^2\big(z,u^1\big)\sin^2\big(u^2\big)+h'^2\big(z,u^2\big)\big)\big(du^2\big)^2,
\\
|d\mathcal{X}|^2
=\cos^2\big(u^2\big)\big(a_1-(a_1-a_0)\sin^2\big(u^1\big)\big)\big(du^1\big)^2\\
\phantom{|d\mathcal{X}|^2=}{}
+\tfrac{1}{2}d\big(\cos^2\big(u^2\big)\big)
d\big(a_0+(a_1-a_0)\sin^2\big(u^1\big)\big)
\\
\phantom{|d\mathcal{X}|^2=}{}
+\big(a_2-\big(a_2-a_0-(a_1-a_0)\sin^2\big(u^1\big)\big)\sin^2\big(u^2\big)\big)\big(du^2\big)^2.
\end{gather*}

Thus the condition $|d\mathcal{X}_z|^2=|d\mathcal{X}|^2$ becomes
\begin{gather*}
f^2\big(z,u^1\big)+\big(a_2-a_0-(a_1-a_0)\sin^2\big(u^1\big)\big)=\mathrm{const}=
\frac{a_2-h'^2(z,u^2)}{\sin^2(u^2)},
\\
f'^2\big(z,u^1\big)+f^2\big(z,u^1\big)g'^2\big(z,u^1\big)=a_1-(a_1-a_0)\sin^2\big(u^1\big),
\end{gather*}
from where we get
\begin{gather}
h\big(z_1,u^2\big):=\int_0^{u^2}\sqrt{a_2-(a_2-z_1a_0)\sin^2(t)}dt,\nonumber\\
f\big(z_1,u^1\big):=\sqrt{(1-z_1)a_0+(a_1-a_0)\sin^2(u^1)},\nonumber\\
g\big(z_1,u^1\big):=\int_0^{u^1}\frac{\sqrt{(1-z_1)a_0a_1+(a_1-a_0)z_1a_0\sin^2(t)}}
{(1-z_1)a_0+(a_1-a_0)\sin^2(t)}dt.\label{eq:fg}
\end{gather}
Note that
\begin{gather}
f\big(0,u^1\big)\cos\big(g\big(0,u^1\big)\big)=\sqrt{a_0}\cos\big(u^1\big),\qquad
f\big(0,u^1\big)\sin\big(g\big(0,u^1\big)\big)=\sqrt{a_1}\sin\big(u^1\big),\label{eq:fg0}
\end{gather}
(we assume simplif\/ications of the form
$\sqrt{a}\sqrt{b}\simeq\sqrt{ab}$ with $\sqrt{\cdot}$ having the
usual def\/inition
$\sqrt{re^{i\theta}}:=\sqrt{r}e^{\frac{i\theta}{2}}$, $r>0$,
$-\pi<\theta\le\pi$, since the possible signs are accounted by
symmetries in the principal planes for quadrics and disappear at
the level of the linear element for their deformations), so~$\mathcal{X}=\mathcal{X}_0$.

The coordinates $\mathbf{x}_0$, $\mathbf{x}_1$, $\mathbf{x}_2$ of
$\mathcal{X}_1$ satisfy (modulo a sign at the second formula)
Peterson's formulae:
\begin{gather}
\sqrt{(\mathbf{x}_0)^2+(\mathbf{x}_1)^2}=\sqrt{a_1-a_0}\cos\big (u^2\big) \sin\big(u^1\big),\nonumber\\
\tan^{-1}\left(\frac{\mathbf{x}_1}{\mathbf{x}_0}\right)=
\frac{\sqrt{a_0}}{\sqrt{a_1-a_0}}\tanh^{-1}\big(\cos\big(u^1\big)\big),\nonumber\\
\mathbf{x}_2=\int_0^{u^2}\sqrt{a_2-(a_2-a_0)\sin^2(t)}dt.\label{eq:pf2}
\end{gather}

More generally
\begin{gather}
h\big(z_1,u^2\big):=\int_0^{u^2}\sqrt{h'^2(t)-z_1\sin^2(t)}dt,\nonumber\\
f\big(z_1,u^1\big):=\sqrt{z_1+f^2(u^1)},\nonumber\\
g\big(z_1,u^1\big):=\int_0^{u^1}\frac{\sqrt{z_1(f'^2(t)+f^2(t)g'^2(t))+f^4(t)g'^2(t)}}
{z_1+f^2(t)}dt\label{eq:bla}
\end{gather}
give Peterson's $1$-dimensional family of deformations~(\ref{eq:surf}) with common conjugate system $(u^1,u^2)$.

\section{Peterson's deformations of higher dimensional quadrics}\label{sec:pet1}

Again we shall discuss only the case of quadrics with center and
having distinct eigenvalues of the quadratic part def\/ining the
quadric, without insisting on totally real cases and deformations
(when the linear elements are real valued).

\begin{remark}
A metric classif\/ication of all (totally real) quadrics in
$\mathbb{C}^{n+1}$ requires the notion of {\it symmetric Jordan}
canonical form of a symmetric complex matrix (see, e.g.~\cite{H}). The symmetric Jordan blocks are:
\begin{gather*}
J_1:=0=0_{1,1}\in\mathbf{M}_1(\mathbb{C}),\qquad
J_2:=f_1f_1^T\in\mathbf{M}_2(\mathbb{C}),\qquad
J_3:=f_1e_3^T+e_3f_1^T\in\mathbf{M}_3(\mathbb{C}),
\\
J_4:=f_1\bar f_2^T+f_2f_2^T+\bar
f_2f_1^T\in\mathbf{M}_4(\mathbb{C}),\qquad J_5:= f_1\bar
f_2^T+f_2e_5^T+e_5f_2^T+\bar
f_2f_1^T\in\mathbf{M}_5(\mathbb{C}),
\\
J_6:= f_1\bar f_2^T+f_2\bar f_3^T+f_3f_3^T+\bar f_3f_2^T+\bar
f_2f_1^T\in\mathbf{M}_6(\mathbb{C}),
\end{gather*}
etc., where $f_j:=\frac{e_{2j-1}-ie_{2j}}{\sqrt{2}}$ are the
standard isotropic vectors (at least the blocks $J_2$, $J_3$ were
known to the classical geometers). Any symmetric complex matrix
can be brought via conjugation with a complex rotation to the
symmetric Jordan canonical form, that is a matrix block
decomposition with blocks of the form $a_jI_p+J_p$; totally real
quadrics are obtained for eigenvalues $a_j$ of the quadratic part
def\/ining the quadric being real or coming in complex conjugate
pairs $a_j$, $\bar a_j$ with subjacent symmetric Jordan blocks of
same dimension $p$.
\end{remark}

Consider the quadric $\sum\limits_{j=0}^n\frac{(x_j)^2}{a_j}=1$, $
a_j\in\mathbb{C}^*$ distinct with parametrization given by the
spherical coordinates on the unit sphere
$\mathbb{S}^n\subset\mathbb{C}^{n+1}$
\[
\mathcal{X}=\sqrt{a_0}\mathbf{C}_0e_0+
\sum_{k=1}^n\sqrt{a_k}\mathbf{C}_k\sin\big(u^k\big)e_k,\qquad
\mathbf{C}_k:=\prod_{j=k+1}^n\cos\big(u^j\big).
\]
The correct generalization of (\ref{eq:surf}) allows us to build
Peterson's deformations of higher dimensional quadrics. With an
eye to the case $n=2$ we make the natural ansatz
\begin{gather}
\mathcal{X}_{\mathbf{z}}=\sum_{k=1}^{n-1}\mathbf{C}_kf_k\big(\mathbf{z},u^k\big)
\big(\cos\big(g_k\big(\mathbf{z},u^k\big)\big)e_{2k-2}
+\sin\big(g_k\big(\mathbf{z},u^k\big)\big)e_{2k-1}\big)+h\big(\mathbf{z},u^n\big)e_{2n-2}\label{eq:surf3}
\end{gather}
with the parameter(s) $\mathbf{z}=(z_1,z_2,\dots )$ to be determined
later.

We have
\begin{gather*}
|d\mathcal{X}_{\mathbf{z}}|^2=\sum_{k=1}^{n-1}\big[\mathbf{C}_k^2\big(f_k'^2\big(\mathbf{z},u^k\big)
+f_k^2\big(\mathbf{z},u^k\big)g_k'^2\big(\mathbf{z},u^k\big)\big)\big(du^k\big)^2
+\tfrac{1}{2}d\big(\mathbf{C}_k^2\big)d\big(f_k^2\big(\mathbf{z},u^k\big)\big)
\\
\phantom{|d\mathcal{X}_{\mathbf{z}}|^2=}{}
+f_k^2(\mathbf{z},u^k)(d\mathbf{C}_k)^2\big]
+h'^2\big(\mathbf{z},u^n\big)\big(du^n\big)^2,\\
|d\mathcal{X}|^2=a_0(d\mathbf{C}_0)^2
+\sum_{k=1}^na_k\big(d\big(\mathbf{C}_k\sin\big(u^k\big)\big)\big)^2.
\end{gather*}
Comparing the coef\/f\/icients of $(du^n)^2$ from
$|d\mathcal{X}_{\mathbf{z}}|^2=|d\mathcal{X}|^2$ we get
\begin{gather*}
\frac{1}{\cos^2(u^n)}\left[\mathbf{C}_1^2\big(f_1^2\big(\mathbf{z},u^1\big)-a_0-(a_1-a_0)\sin^2\big(u^1\big)\big)+
\sum_{k=2}^{n-1}\mathbf{C}_k^2\big(f_k^2\big(\mathbf{z},u^k\big)-a_k\sin^2\big(u^k\big)\big)\right]
\\
\qquad {}=\mathrm{const}=\frac{a_n\cos^2(u^n)-h'^2(\mathbf{z},u^n)}{\sin^2(u^n)}
\end{gather*}
from where we get with $z_0:=1$:
\begin{gather*}
f_k^2\big(z_{k-1},z_k,u^k\big):=(z_{k-1}-z_k)a_0+(a_k-z_{k-1}a_0)\sin^2\big(u^k\big),\qquad
k=1,\dots ,n-1,
\\
h'^2\big(z_{n-1},u^n\big):=a_n-(a_n-z_{n-1}a_0)\sin^2\big(u^n\big).
\end{gather*}
Now we have
\begin{gather*}
(d\mathbf{C}_0)^2=\sum_{k=1}^{n-1}\big[z_{k-1}(d\mathbf{C}_{k-1})^2-z_k(d\mathbf{C}_k)^2\big]
+z_{n-1}(d\mathbf{C}_{n-1})^2=\sum_{k=1}^{n-1}\!\big[z_{k-1} \big(\mathbf{C}_k^2\sin^2\big(u^k\big)\big(du^k\big)^2
\\
\phantom{(d\mathbf{C}_0)^2=}{}
-\tfrac{1}{2}d\big(\mathbf{C}_k^2\big)d\big(\sin^2(u^k)\big)+\cos^2\big(u^k\big)\big(d\mathbf{C}_k\big)^2\big)
-z_k\big(d\mathbf{C}_k\big)^2\big] +z_{n-1}\big(d\mathbf{C}_{n-1}\big)^2,
\\
\big(d\big(\mathbf{C}_k\sin\big(u^k\big)\big)\big)^2=\mathbf{C}_k^2\cos^2\big(u^k\big)\big(du^k\big)^2
+\tfrac{1}{2}d\big(\mathbf{C}_k^2\big)d\big(\sin^2\big(u^k\big)\big)+\sin^2\big(u^k\big)\big(d\mathbf{C}_k\big)^2,
\end{gather*}
so
\begin{gather*}
|d\mathcal{X}|^2=\sum_{k=1}^{n-1}\big[\mathbf{C}_k^2\big(a_k-
(a_k-z_{k-1}a_0)\sin^2\big(u^k\big)\big)\big(du^k\big)^2
+\tfrac{1}{2}(a_k-z_{k-1}a_0)d\big(\mathbf{C}_k^2\big)d\big(\sin^2\big(u^k\big)\big)
\\
\phantom{|d\mathcal{X}|^2=}{}
+\big((z_{k-1}-z_k)a_0+
(a_k-z_{k-1}a_0)\sin^2\big(u^k\big)\big)\big(d\mathbf{C}_k\big)^2\big]\\
\phantom{|d\mathcal{X}|^2=}{}
+\big(a_n-(a_n-z_{n-1}a_0)\sin^2\big(u^n\big)\big)\big(du^n\big)^2,
\\
0=|d\mathcal{X}_{\mathbf{z}}|^2-|d\mathcal{X}|^2=
\sum_{k=1}^{n-1}\mathbf{C}_k^2\big(f_k'^2\big(\mathbf{z},u^k\big)
+f_k^2\big(\mathbf{z},u^k\big)g_k'^2\big(\mathbf{z},u^k\big)\\
\phantom{0=}{}-a_k+
(a_k-z_{k-1}a_0)\sin^2\big(u^k\big)\big)\big(du^k\big)^2,
\end{gather*}
so we f\/inally get (\ref{eq:fgn}).

For $z_1=z_2=\dots =z_{n-1}=0$ we get $g_2=\dots =g_{n-1}=0$ and using
(\ref{eq:fg0}) we get $\mathcal{X}=\mathcal{X}_{\mathbf{0}}$ with
$\mathbb{C}^{n+1}\hookrightarrow\mathbb{C}^{2n-1}$ as
$(x_0,x_1,\dots ,x_n)\mapsto(x_0,x_1,x_2,0,x_3,0,\dots ,x_{n-1},0,x_n)$.

For $z_1=z_2=\dots =z_{n-1}=1$ we get
$\mathcal{X}_{\mathbf{1}}=(\mathbf{x}_0,\dots ,\mathbf{x}_{2n-2})$
given by Peterson's formulae (\ref{eq:pfn}).

More generally and with $z_0:=0$
\begin{gather}
f_k\big(z_{k-1},z_k,u^k\big):=\sqrt{z_k+f_k^2(u^k)-z_{k-1}\cos^2(u^k)},\qquad k=1,\dots ,n-1,\nonumber\\
g_k\big(z_{k-1},z_k,u^k\big):=\int_0^{u^k}\frac{\sqrt{f_k'^2(t)+f_k^2(t)g_k'^2(t)-
(f_k'^2(z_{k-1},z_k,t)+z_{k-1}\sin^2(t))}}
{f_k(z_{k-1},z_k,t)}dt,\nonumber\\
h\big(z_{n-1},u^n\big):=\int_0^{u^n}\sqrt{h'^2(t)-z_{n-1}\sin^2(t)}dt \label{eq:fgr}
\end{gather}
give an $(n-1)$-dimensional family of deformations
(\ref{eq:surf3}); for $g_k(u^k)=0$, $k=2,\dots ,n-1$ we have
$\mathcal{X}_{\mathbf{0}}\subset\mathbb{C}^{n+1}$.

\section{The common conjugate system\\ and non-degenerate joined second fundamental forms}\label{sec:pet2}

The fact that $(u^1,\dots ,u^n)$ is a conjugate system on
$\mathcal{X}_{\mathbf{0}}$ is clear since we have
\[
\partial_{u^k}\partial_{u^j}\mathcal{X}_\mathbf{0}=
-\tan\big(u^j\big)\partial_{u^k}\mathcal{X}_\mathbf{0},\qquad 1\le k<j\le n.
\]
With the normal f\/ield
\[
\hat
N_{\mathbf{0}}:=(\sqrt{a_0})^{-1}\mathbf{C}_0e_0+
\sum_{k=1}^n(\sqrt{a_k})^{-1}\mathbf{C}_k\sin\big(u^k\big)e_k
\]
we have $\hat
N_{\mathbf{0}}^Td^2\mathcal{X}_{\mathbf{0}}=-\sum\limits_{k=1}^n\mathbf{C}_k^2(du^k)^2$.
To see that $(u^1,\dots ,u^n)$ is a conjugate system on
\begin{gather*}
\mathcal{X}=(x_0,\dots ,x_{2n-2}):=\sum_{k=1}^{n-1}\mathbf{C}_kf_k\big(u^k\big)\big(\cos\big(g_k\big(u^k\big)\big)e_{2k-2}
+\sin\big(g_k\big(u^k\big)\big)e_{2k-1}\big)\\
\phantom{\mathcal{X}=(x_0,\dots ,x_{2n-2}):=}{} +h\big(u^n\big)e_{2n-2}
\end{gather*}
we have again
$\partial_{u^k}\partial_{u^j}\mathcal{X}=-\tan(u^j)\partial_{u^k}\mathcal{X}$, $1\le k<j\le n$; again the $n-1$ f\/ields
\begin{gather*}
\partial_{u^1}\mathcal{X}|_{u^1,u^2,\dots ,\widehat{u^k},\dots ,u^n=\mathrm{const}},\qquad
\partial_{u^2}\mathcal{X}|_{u^1,u^2,\dots ,\widehat{u^k},\dots ,u^n=\mathrm{const}},\qquad
\dots ,\\
\widehat{\partial_{u^k}\mathcal{X}}|_{u^1,u^2,\dots ,\widehat{u^k},\dots ,u^n=\mathrm{const}},\qquad
\dots ,
\qquad
\partial_{u^n}\mathcal{X}|_{u^1,u^2,\dots ,\widehat{u^k},\dots ,u^n=\mathrm{const}},\qquad
k=1,\dots ,n
\end{gather*}
generate ruled $n$-dimensional developables in $\mathbb{C}^{2n-1}$
because the only term producing shape is
$\partial_{u^k}\partial_{u^k}\mathcal{X}$.

For the non-degenerate joined second fundamental forms property we
have
\begin{gather*}
u^n=h^{-1}(x_{2n-2}),\qquad h'\big(u^n\big)du^n=dx_{2n-2},\qquad
u^k=g_k^{-1}\left(\tan^{-1}\left(\frac{x_{2k-1}}{x_{2k-2}}\right)\right),
\\
\mathbf{C}_k^2f_k^2(u^k)g'_k(u^k)du^k=x_{2k-2}dx_{2k-1}-x_{2k-1}dx_{2k-2},
\qquad k=1,\dots ,n-1
\end{gather*}
and $\mathcal{X}$ is given implicitly by the zeroes of the
functionally independent
\[
F_k:=(x_{2k-2})^2+(x_{2k-1})^2-\mathbf{C}_k^2f_k^2\big(u^k\big),\qquad
k=1,\dots ,n-1.
\]
We have the natural linearly independent normal f\/ields
\begin{gather*}
N_k:=\nabla F_k=2(x_{2k-2}e_{2k-2}+x_{2k-1}e_{2k-1})
-\frac{2f_k'(u^k)(-x_{2k-1}e_{2k-2}+x_{2k-2}e_{2k-1})}{f_k(u^k)g'_k(u^k)}
\\
\phantom{N_k:=}{} +2\mathbf{C}_k^2f_k^2(u^k)\left[\sum_{j=k+1}^{n-1}\frac{\tan(u^j)
(-x_{2j-1}e_{2j-2}+x_{2j-2}e_{2j-1})}
{\mathbf{C}_j^2f_j^2(u^j)g'_j(u^j)}+\frac{\tan(u^n)e_{2n-2}}{h'(u^n)}\right],\\
\phantom{N_k:=}{} k=1,\dots ,n-1,
\end{gather*} and
\begin{gather*}
\partial_{u^l}\partial_{u^l}\mathcal{X}=-\sum_{j=1}^{l-1}(x_{2j-2}e_{2j-2}+x_{2j-1}e_{2j-1})
+\left(\frac{f_l''(u^l)}{f_l(u^l)}-g_l'^2\big(u^l\big)\right)(x_{2l-2}e_{2l-2}+x_{2l-1}e_{2l-1})
\\
\phantom{\partial_{u^l}\partial_{u^l}\mathcal{X}=}{} +g'_l\big(u^l\big)\left(\frac{2f'_l(u^l)}{f_l(u^l)}+\frac{g_l''(u^l)}{g_l'(u^l)}\right)
(-x_{2l-1}e_{2l-2}+x_{2l-2}e_{2l-1}),\qquad l=1,\dots ,n-1,
\\
\partial_{u^n}\partial_{u^n}\mathcal{X}=-\sum_{l=1}^{n-1}(x_{2l-2}e_{2l-2}+x_{2l-1}e_{2l-1})
+h''\big(u^n\big)e_{2n-2},
\end{gather*}
and the second fundamental form
\begin{gather*}
N_k^Td^2\mathcal{X}=2\mathbf{C}_k^2f_k^2\Bigg[\left(\frac{f_k''(u^k)}{f_k(u^k)}-g_k'^2\big(u^k\big)
-\frac{f_k'(u^k)}{f_k(u^k)}\left(\frac{2f'_k(u^k)}{f_k(u^k)}+
\frac{g_k''(u^k)}{g_k'(u^k)}\right)\right)\big(du^k\big)^2
\\
\phantom{N_k^Td^2\mathcal{X}=}{}
+\sum_{l=k+1}^{n-1}\left(\tan\big(u^l\big)\left(\frac{2f'_l(u^l)}{f_l(u^l)}
+\frac{g_l''(u^l)}{g_l'(u^l)}\right)-1\right)\big(du^l\big)^2\\
\phantom{N_k^Td^2\mathcal{X}=}{}
+\left(\frac{\tan(u^n)h''(u^n)}{h'(u^n)}-1\right)\big(du^n\big)^2\Bigg],\qquad k=1,\dots ,n-1.
\end{gather*}

For Peterson's deformations of higher dimensional quadrics we have
\begin{gather*}
N_k^Td^2\mathcal{X}=
-2a_0\mathbf{C}_k^2f_k^2\left(\frac{a_kz_{k-1}(du^k)^2}{g_k'^2(u^k)f_k^4(u^k)}
+\sum_{l=k+1}^{n-1}\frac{a_l(z_{l-1}-z_l)(du^l)^2}{g_l'^2(u^l)f_l^4(u^l)}
+\frac{a_n(du^n)^2}{a_0h'^2(u^n)}\right).
\end{gather*}
It is now enough to check the open non-degenerate joined second
fundamental forms property only for $\mathbf{z}=(1,1,\dots ,1)$. Thus
with $\delta:=\frac{a_n}{a_0\sin^2(u^n)+a_n\cos^2(u^n)}$ we need
\begin{gather*}
0\neq\begin{vmatrix}
\mathbf{C}_1&\mathbf{C}_2&\mathbf{C}_3&\dots &\mathbf{C}_{n-1}&\delta^{-1}\mathbf{C}_n\\
\frac{a_1}{a_1-a_0}&0&0&\dots &0&\sin^2\big(u^1\big)\\
0&\frac{a_2}{a_2-a_0}&0&\dots &0&\sin^2\big(u^2\big)\\
0&0&\frac{a_3}{a_3-a_0}&\dots &0&\sin^2\big(u^3\big)\\
\vdots&\vdots&\vdots&\cdots&\vdots&\vdots\\
0&0&0&\dots &\frac{a_{n-1}}{a_{n-1}-a_0}&\sin^2\big(u^{n-1}\big)\end{vmatrix}
\end{gather*}
almost everywhere, which is straightforward.

\begin{remark}
Note that a-priori $\mathcal{X}_{\mathbf{1}}$ comes close to lie
in a degenerate deformation of $\mathbb{C}^{n+1}$ in
$\mathbb{C}^{2n-1}$:  $\hat
N_{\mathbf{0}}^Td^2\mathcal{X}_{\mathbf{0}}-
\big(\sum\limits_{k=1}^{n-1}\frac{1}{2a_k}N_k\big)^Td^2\mathcal{X}_{\mathbf{1}}$
depends only on $(du^n)^2$ and this is as closest to $0$ as we can
get.
\end{remark}

\section{Conjugate systems}\label{sec:pet3}

Consider the complexif\/ied Euclidean space
\[
\big(\mathbb{C}^n,\langle \cdot ,\cdot \rangle),\qquad \langle x,y\rangle :=x^Ty,\qquad |x|^2:=x^Tx,\qquad
x,y\in\mathbb{C}^n
\]
with standard basis $\{e_j\}_{j=1,\dots ,n}$, $
e_j^Te_k=\delta_{jk}.$

Isotropic (null) vectors are those vectors $v$ of length $0$  $(|v|^2=0)$; since most vectors are not isotropic we shall call a
vector simply vector and we shall only emphasize isotropic when
the vector is assumed to be isotropic. The same denomination will
apply in other settings: for example we call quadric a
non-degenerate quadric (a quadric projectively equivalent to the
complex unit sphere).

For $n\ge 3$ consider the $n$-dimensional sub-manifold
\[
 x=x\big(u^1,u^2,\dots ,u^n\big)\subset\mathbb{C}^{n+p}, \qquad
du^1\wedge du^2\wedge\cdots \wedge du^n\neq 0
\]
such that the tangent
space at any point of $x$ is not isotropic (the scalar product
induced on it by the Euclidean one on $\mathbb{C}^{n+p}$ is not
degenerate; this assures the existence of orthonormal normal
frames). We shall always have Latin indices
$j,k,l,m,p,q\in\{1,\dots ,n\}$, Greek ones
$\alpha,\beta,\gamma\in\{n+1,\dots ,n+p\}$ and mute summation for upper and
lower indices when clear from the context; also we shall preserve
the classical notation $d^2$ for the tensorial (symmetric) second
derivative. We have the normal frame $N:=[N_{n+1} \ \dots \
N_{n+p}],\ N^TN=I_p$, the f\/irst $|dx|^2=g_{jk}du^j\odot du^k$ and
second $d^2x^TN=[h_{jk}^{n+1}du^j\odot du^k\ \dots \
h_{jk}^{n+p}du^j\odot du^k]$ fundamental forms, the Christof\/fel
symbols
$\Gamma_{jk}^l=\frac{g^{lm}}{2}[\partial_{u^k}g_{jm}+\partial_{u^j}g_{km}-\partial_{u^m}g_{jk}]$,
the Riemann curvature
$R_{jmkl}=g_{mp}R^p_{jkl}=g_{mp}[\partial_{u^l}\Gamma_{jk}^p-\partial_{u^k}\Gamma_{jl}^p+\Gamma_{jk}^q\Gamma_{ql}^p
-\Gamma_{jl}^q\Gamma_{qk}^p]$ tensor, the normal connection
$N^TdN=\{n^{\alpha}_{\beta j}du^j\}_{\alpha,\beta=n+1,\dots ,n+p}$, $n^{\alpha}_{\beta j}=-n^{\beta}_{\alpha j}$ and the curvature
$r_{\alpha jk}^{\beta}=\partial_{u^k}n_{\alpha
j}^{\beta}-\partial_{u^j}n_{\alpha k}^{\beta}+n_{\alpha
j}^{\gamma}n_{\gamma k}^{\beta}-n_{\alpha k}^{\gamma}n_{\gamma
j}^{\beta}$ tensor of the normal bundle.

We have the {\it Gau\ss--Weingarten} (GW) equations
\[
\partial_{u^k}\partial_{u^j}x=\Gamma_{jk}^l\partial_{u^l}x+h_{jk}^{\alpha}N_{\alpha},\qquad
\partial_{u^j}N_{\alpha}=
-h_{jk}^{\alpha}g^{kl}\partial_{u^l}x+n^{\beta}_{j\alpha}N_{\beta}
\]
and their integrability conditions
$\partial_{u^l}(\partial_{u^k}\partial_{u^j}x)=\partial_{u^k}(\partial_{u^l}\partial_{u^j}x)$,
$\partial_{u^k}(\partial_{u^j}N_{\alpha})=\partial_{u^j}(\partial_{u^k}N_{\alpha})$, from
where one obtains by taking the tangential and normal components
(using
$-\partial_{u^l}g^{jk}=g^{jm}\Gamma_{ml}^k+g^{km}\Gamma_{ml}^j$
and the GW equations themselves) the {\it
Gau\ss--Codazzi--Mainardi$($--Peterson$)$--Ricci} \mbox{(G-CMP-R)} equations
\begin{gather*}
R_{jmkl}=\sum_{\alpha}\big(h_{jk}^{\alpha}h_{lm}^{\alpha}-h_{jl}^{\alpha}h_{km}^{\alpha}\big),\\
\partial_{u^l}h_{jk}^{\alpha}-\partial_{u^k}h_{jl}^{\alpha}
+\Gamma_{jk}^mh_{ml}^{\alpha}-\Gamma_{jl}^mh_{mk}^{\alpha}
+h_{jk}^{\beta}n_{\beta l}^{\alpha}-h_{jl}^{\beta}n_{\beta
k}^{\alpha}=0,
\\
r_{\alpha
jk}^{\beta}=h_{jl}^{\alpha}g^{lm}h_{mk}^{\beta}-h_{kl}^{\alpha}g^{lm}h_{mj}^{\beta}.
\end{gather*}
If we have conjugate system
$h_{jk}^{\alpha}=:\delta_{jk}h_j^{\alpha}$, then the above
equations become:
\begin{gather}
R_{jkjk}=-R_{jkkj}=\sum_{\alpha}h_j^{\alpha}h_k^{\alpha},\qquad
\partial_{u^k}h_j^{\alpha}=\Gamma_{jk}^jh_j^{\alpha}-\Gamma_{jj}^kh_k^{\alpha}-h_j^{\beta}\eta_{\beta
k}^{\alpha},\quad j\neq k,\nonumber\\
R_{jklm}=0\quad \mathrm{otherwise},\nonumber\\
\Gamma_{jk}^lh_l^{\alpha}=\Gamma_{jl}^kh_k^{\alpha},\quad  j,k,l\
\mathrm{distinct},\qquad r_{\alpha
jk}^{\beta}=(h_j^{\alpha}h_k^{\beta}-h_j^{\beta}h_k^{\alpha})g^{jk}.\label{eq:riem}
\end{gather}
In particular for lines of curvature parametrization
($g_{jk}=\delta_{jk}g_{jk}$) we have f\/lat normal bundle, so one
can choose up to multiplication on the right by a constant matrix
$\in\mathbf{O}_p(\mathbb{C})$ normal frame~$N$ with zero normal
connection $N^TdN=0$.

This constitutes a dif\/ferential system in the $np$ unknowns
$h_j^{\alpha}$ and the yet to be determined coef\/f\/icients
$\eta_{\beta k}^{\alpha}$; according to Cartan's exterior
dif\/ferential systems in involution tools in order to study
deformations of $n$-dimensional sub-manifolds of
$\mathbb{C}^{n+p}$ in conjugate system para\-me\-terization one must
iteratively apply compatibility conditions (commuting of mixed
derivatives) to the equations of this system and their
algebraic-dif\/ferential consequences, introducing new variables as
necessary and assuming only identities obtained at previous
iterations and general identities for the Riemann curvature tensor
(symmetries and Bianchi identities):
\begin{gather*}
R_{jklm}=-R_{kjlm}=-R_{jkml}=R_{lmjk},\qquad
R_{jklm}+R_{jlmk}+R_{jmkl}=0,\\ R_{jklm;q}+R_{jkmq;l}+R_{jkql;m}=0,
\\
R_{jklm;q}:=\partial_{u^q}R_{jklm}-\Gamma_{qj}^rR_{rklm}-\Gamma_{qk}^rR_{jrlm}
-\Gamma_{ql}^rR_{jkrm}-\Gamma_{qm}^rR_{jklr}
\end{gather*}
until no further conditions appear from compatibility conditions.
However one cannot use in full the Cartan's exterior dif\/ferential
forms and moving frames tools (see, e.g.~\cite{B}), since
they are best suited for arbitrary (orthonormal) tangential frames
and orthonormal normal ones and their corresponding change of
frames; thus one loses the advantage of special coordinates suited
to our particular problem.

In our case we only obtain
\begin{gather}
\partial_{u^l}R_{jkjk}=\big(\Gamma_{jl}^j+\Gamma_{kl}^k\big)R_{jkjk}-\Gamma_{kk}^lR_{jljl}-\Gamma_{jj}^lR_{klkl},\quad
j,k,l\ \ \mathrm{distinct},\nonumber\\
\Gamma_{lk}^mR_{jmjm}-\Gamma_{mk}^lR_{jljl}=0,\quad  j,k,l,m \ \
\mathrm{distinct}.\label{eq:riem1}
\end{gather}

\begin{remark}
Dif\/ferentiating the f\/irst equations of (\ref{eq:riem}) with
respect to $u^l$, $l\neq j,k$ and using~(\ref{eq:riem}) itself we
obtain
\begin{gather*}
\partial_{u^l}R_{jkjk}=\sum_{\alpha}
\big(\partial_{u^l}h_j^{\alpha}h_k^{\alpha}+h_j^{\alpha}\partial_{u^l}h_k^{\alpha}\big)\\
\phantom{\partial_{u^l}R_{jkjk}}{}
=\sum_{\alpha}\big[\big(\Gamma_{jl}^jh_j^{\alpha}-\Gamma_{jj}^lh_l^{\alpha}-h_j^{\beta}\eta_{\beta
l}^{\alpha}\big)h_k^{\alpha}+
h_j^{\alpha}\big(\Gamma_{kl}^kh_k^{\alpha}-\Gamma_{kk}^lh_l^{\alpha}-h_k^{\beta}\eta_{\beta
l}^{\alpha}\big)\big]\\
\phantom{\partial_{u^l}R_{jkjk}}{}
=\big(\Gamma_{jl}^j+\Gamma_{kl}^k\big)R_{jkjk}-\Gamma_{kk}^lR_{jljl}-\Gamma_{jj}^lR_{klkl},
\end{gather*}
that is the f\/irst equations of~(\ref{eq:riem1}), so the covariant
derivative of the Gau\ss\ equations become, via the G-CMP
equations, the Bianchi second identity (the second equations of~(\ref{eq:riem1}) being consequence of the CMP equations is obvious; see also~\cite{B}).
\end{remark}

\section{The non-degenerate joined second fundamental forms\\ assumption}\label{sec:pet4}

With an eye towards our interests (deformations in
$\mathbb{C}^{2n-1}$ of quadrics in $\mathbb{C}^{n+1}$ and with
common conjugate system) we make the genericity assumption of
non-degenerate joined second fundamental forms of~$x_0$,~$x$: with
$d^2x_0^TN_0=:h_j^0(du^j)^2$ being the second fundamental form of
the quadric $x_0\subset\mathbb{C}^{n+1}$ whose deformation
$x\subset\mathbb{C}^{2n-1}$ is (that is $|dx_0|^2=|dx|^2$) the
vectors $h_j:=[ih_j^0 \ h_j^{n+1} \ \dots \ h_j^{2n-1}]^T$ are
linearly independent. From the Gau\ss\ equations we obtain
$h_j^0h_k^0=R_{jkjk}=\sum_{\alpha}h_j^{\alpha}h_k^{\alpha}, \  j\neq
k \ \Leftrightarrow \ h_j^Th_k=\delta_{jk}|h_j|^2$; thus the vectors
$h_j\subset\mathbb{C}^n$ are further orthogonal, which prevents
them from being isotropic (should one of them be isotropic, by a~rotation of $\mathbb{C}^n$ one can make it $f_1$ and after
subtracting suitable multiples of $f_1$ from the remaining ones by
another rotation of $\mathbb{C}^n$ the remaining ones linear
combinations of $e_3,\dots ,e_n$, so we would have $n-1$ linearly
independent orthogonal vectors in $\mathbb{C}^{n-2}$, a
contradiction), so $\mathbf{a}_j:=|h_j|\neq 0$,
$h_j=:\mathbf{a}_jv_j$, $R:=[v_1\ \dots \
v_n]\subset\mathbf{O}_n(\mathbb{C})$.

Thus with $\eta_{0j}^{\alpha}=-\eta_{\alpha j}^0:=0$, $(\eta_{\beta
j}^{\alpha})_{\alpha,\beta=0,n+1,\dots ,2n-1}=:\Upsilon_j=-\Upsilon_j^T$
we have reduced the problem to f\/inding $R=[v_1\ \dots \
v_n]\subset\mathbf{O}_n(\mathbb{C})$,
$\mathbf{a}_j\subset\mathbb{C}^*$,
$\Upsilon_j\subset\mathbf{M}_n(\mathbb{C})$,
$\Upsilon_j=-\Upsilon_j^T$, $\Upsilon_je_1=0$ satisfying the
dif\/ferential system
\begin{gather}
\partial_{u^k}\log\mathbf{a}_j=\Gamma_{jk}^j,\qquad
\partial_{u^k}v_j=-\Gamma_{jj}^k\frac{\mathbf{a}_k}{\mathbf{a}_j}v_k-\Upsilon_kv_j,\nonumber\\
\partial_{u^k}\Upsilon_j-\partial_{u^j}\Upsilon_k-[\Upsilon_j,\Upsilon_k]
=-g^{jk}\mathbf{a}_j\mathbf{a}_k\big(I_n-e_1e_1^T\big)\big(v_jv_k^T-v_kv_j^T\big)\big(I_n-e_1e_1^T\big),\quad
j\neq k,\nonumber\\
\sum_j\frac{(h_j^0)^2}{\mathbf{a}_j^2}+1=0\label{eq:syst}
\end{gather} derived from the CMP-R
equations and $\mathbf{a}_jv_j^1=ih_j^0$, $\sum_j(v_j^1)^2=1$ and
with the linear element further satisfying the condition
\begin{gather}
\Gamma_{jk}^l=0,\quad  j,k,l\ \ \mathrm{distinct}\label{eq:jkl}
\end{gather}
derived from the CMP equations.

First we shall investigate the consequences of (\ref{eq:jkl}), via
the properties of the Riemann curvature tensor, on the other
Christof\/fel symbols. For $j$, $k$, $l$ distinct we have
$0=g^{pm}R_{jmkl}=R^p_{jkl}=\partial_{u^l}\Gamma_{jk}^p-\partial_{u^k}\Gamma_{jl}^p+\Gamma_{jk}^q\Gamma_{ql}^p
-\Gamma_{jl}^q\Gamma_{qk}^p$; thus for $p=k$ we obtain
\begin{gather}
\partial_{u^l}\Gamma_{kj}^k=\Gamma_{kl}^k\Gamma_{lj}^l+\Gamma_{kj}^k\Gamma_{jl}^j-\Gamma_{kl}^k\Gamma_{kj}^k,\quad
j,k,l\ \ \mathrm{distinct}.\label{eq:jkla}
\end{gather}
We also have $R^p_{jjl}=g^{pm}R_{jmjl}=g^{pl}R_{jljl}$, $j\neq l$,
so
\begin{gather}
g^{pl}R_{jljl}=\partial_{u^l}\Gamma_{jj}^p+\Gamma_{jj}^p\big(\Gamma_{pl}^p-\Gamma_{jl}^j\big)
+\Gamma_{jj}^l\Gamma_{ll}^p,\quad j,l,p\ \
\mathrm{distinct},\nonumber\\
g^{jl}R_{jljl}=\partial_{u^l}\Gamma_{jj}^j-\partial_{u^j}\Gamma_{jl}^j
+\Gamma_{jj}^l\Gamma_{ll}^j
-\Gamma_{lj}^l\Gamma_{jl}^j,\nonumber\\
g^{ll}R_{jljl}=\partial_{u^l}\Gamma_{jj}^l-\partial_{u^j}\Gamma_{lj}^l
+\Gamma_{jj}^q\Gamma_{lq}^l-\Gamma_{jl}^j\Gamma_{jj}^l
-\big(\Gamma_{lj}^l\big)^2,\quad j\neq l.\label{eq:jjl}
\end{gather}
Conversely, (\ref{eq:jkl}) and (\ref{eq:jkla}) imply $R_{jklm}=0$
for three of $j$, $k$, $l$, $m$ distinct.

\begin{remark}
Note that (\ref{eq:jkl}) are valid for orthogonal coordinates, so
conjugate systems with the property (\ref{eq:jkl}) are a natural
projective generalization of lines of curvature on $n$-dimensional
sub-manifolds $x\subset\mathbb{C}^{n+p}$ (to see this f\/irst
$\partial_{u^k}\partial_{u^j}x=\Gamma_{jk}^j\partial_{u^j}x+\Gamma_{kj}^k\partial_{u^k}x$,
$j\neq k$ is af\/f\/ine invariant (thus $\Gamma_{jk}^j$, $\Gamma_{kj}^k$
are also af\/f\/ine invariants) and further
$\partial_{u^k}\partial_{u^j}\frac{x}{\rho}=(\Gamma_{jk}^j-\partial_{u^k}\log\rho)\partial_{u^j}\frac{x}{\rho}
+(\Gamma_{kj}^k-\partial_{u^j}\log\rho)\partial_{u^k}\frac{x}{\rho}$, $
j\neq k$ for $\rho\subset\mathbb{C}^*$ with
$\partial_{u^k}\partial_{u^j}\rho=\Gamma_{jk}^j\partial_{u^j}\rho+\Gamma_{kj}^k\partial_{u^k}\rho$, $j\neq k$).
\end{remark}

Imposing the compatibility conditions
\[
\partial_{u^l}(\partial_{u^k}\log\mathbf{a}_j)=\partial_{u^k}(\partial_{u^l}\log\mathbf{a}_j),\qquad
\partial_{u^l}(\partial_{u^k}v_j)=\partial_{u^k}(\partial_{u^l}v_j),\quad  j,k,l \ \ {\rm  distinct}
\]
we obtain
\begin{gather}
\partial_{u^l}\Gamma_{jk}^j-\partial_{u^k}\Gamma_{jl}^j=0=
\partial_{u^l}\Gamma_{jj}^k+\Gamma_{jj}^k\big(\Gamma_{kl}^k
-\Gamma_{jl}^j\big)+\Gamma_{jj}^l\Gamma_{ll}^k-g^{kl}h_j^0h_l^0,\quad
j,k,l\ \ \mathrm{distinct},\label{eq:com}
\end{gather}
which are consequences of (\ref{eq:jkla}) and the f\/irst equations
of (\ref{eq:jjl}).

From the f\/irst equations of (\ref{eq:syst}) we get by integration
a precise determination of $\mathbf{a}_j$ up to multiplication by
a function of $u^j$; thus the $\mathbf{a}_j$ part of the solution
depends on at most $n$ functions of one variable (with the last
equation of (\ref{eq:syst}) also to be taken into consideration);
we shall see later that the remaining part of the dif\/ferential
system involves the normal bundle and its indeterminacy, so it
will not produce a bigger space of solutions (no functional
information is allowed in the normal bundle). From the CMP
equations of~$x_0$ and the second equations of~(\ref{eq:com}) we
obtain with $\gamma_{jk}:=\Gamma_{jj}^k\frac{h_k^0}{h_j^0}$, $j\neq
k$:
\begin{gather}
\partial_{u^l}\gamma_{jk}=\gamma_{jk}\gamma_{jl}-\gamma_{jl}\gamma_{lk}
-\gamma_{jk}\gamma_{kl}+g^{kl}h_k^0h_l^0,\quad j,k,l\ \
\mathrm{distinct}.\label{eq:gaj}
\end{gather}
From the f\/irst equations of (\ref{eq:syst}) and the CMP equations
of $x_0$ we obtain with $\mathbf{b}_j:=\frac{h_j^0}{\mathbf{a}_j}$
\begin{gather}
\partial_{u^k}\log\mathbf{b}_j=-\gamma_{jk},\quad j\neq k,\label{eq:hoj}
\end{gather}
so dif\/ferentiating the last equation of (\ref{eq:syst}) we obtain
\begin{gather}
\partial_{u^j}\log\mathbf{b}_j=\mathbf{b}_j^{-2}\sum_{l\neq
j}\mathbf{b}_l^2\gamma_{lj}.\label{eq:ajj}
\end{gather}
This assures that
\begin{gather}
\Upsilon_k:=-\partial_{u^k}RR^T-\sum_{j\neq
k}\gamma_{jk}\frac{\mathbf{b}_j}{\mathbf{b}_k}\big(v_kv_j^T-v_jv_k^T\big)\label{eq:ups}
\end{gather}
satisf\/ies $e_1^T\Upsilon_k=0$. Thus we have reduced the problem to
f\/inding $\mathbf{b}_j$ satisfying (\ref{eq:hoj}), (\ref{eq:ajj})
(in this case $\sum_j\mathbf{b}_j^2+1=0$ becomes a prime integral
of (\ref{eq:hoj}), (\ref{eq:ajj}) and removes a constant from the
space of solutions) and then completing $v_j^1=-i\mathbf{b}_j$ to
$R=[v_1\dots v_n]\subset\mathbf{O}_n(\mathbb{C})$ in an arbitrary
manner (that is undetermined up to multiplication on the left with
$R'\subset\mathbf{O}_n(\mathbb{C})$, $R'e_1=e_1$); with the second
fundamental form of~$x$ found one f\/inds~$x$ by the integration of
a Ricatti equation and quadratures (the Gau\ss--Bonnet(--Peterson)
theorem). $\Upsilon_j$ given by~(\ref{eq:ups}) will satisfy
\begin{gather}
\partial_{u^k}\Upsilon_j-\partial_{u^j}\Upsilon_k-[\Upsilon_j,\Upsilon_k]
=-\frac{g^{jk}h_j^0h_k^0}{\mathbf{b}_j\mathbf{b}_k}\big(I_n-e_1e_1^T\big)\big(v_jv_k^T-v_kv_j^T\big)
\big(I_n-e_1e_1^T\big),\!\quad j\neq k.\!\!\!\label{eq:upsj}
\end{gather}
Imposing the compatibility condition
$\partial_{u^k}(\partial_{u^j}\log\mathbf{b}_j)=\partial_{u^j}(\partial_{u^k}\log\mathbf{b}_j)$, $j\neq k$ on (\ref{eq:hoj}), (\ref{eq:ajj}) we obtain
\begin{gather}
\partial_{u^k}\left(\gamma_{kj}\frac{\mathbf{b}_k}{\mathbf{b}_j}\right)
+\partial_{u^j}\left(\gamma_{jk}\frac{\mathbf{b}_j}{\mathbf{b}_k}\right)-\sum_{l\neq
j,k}\big(\gamma_{lj}\gamma_{lk}-g^{jk}h_j^0h_k^0\big)\frac{\mathbf{b}_l^2}{\mathbf{b}_j\mathbf{b}_k}=0,\quad
j\neq k.\label{eq:sym}
\end{gather}
Now (\ref{eq:upsj}) becomes
\begin{gather*}
\sum_{m\neq
j}\partial_{u^k}\left(\gamma_{mj}\frac{\mathbf{b}_m}{\mathbf{b}_j}\right)\big(v_jv_m^T-v_mv_j^T\big)
-\sum_{l\neq
k}\partial_{u^j}\left(\gamma_{lk}\frac{\mathbf{b}_l}{\mathbf{b}_k}\right)\big(v_kv_l^T-v_lv_k^T\big)
\\
\qquad{}
+\sum_{m\neq j,\ l\neq
k}\gamma_{mj}\gamma_{lk}\frac{\mathbf{b}_m\mathbf{b}_l}{\mathbf{b}_j\mathbf{b}_k}
\big[\delta_{mk}\big(v_jv_l^T-v_lv_j^T\big)+\delta_{ml}\big(v_kv_j^T-v_jv_k^T\big)
-\delta_{jl}\big(v_kv_m^T-v_mv_k^T\big)\big]
\\
\qquad{}
=\frac{g^{jk}h_j^0h_k^0}{\mathbf{b}_j\mathbf{b}_k}\big(I_n-e_1e_1^T\big)\big(v_jv_k^T-v_kv_j^T\big)
\big(I_n-e_1e_1^T\big),\quad j\neq k,
\end{gather*}
which boils down to (\ref{eq:sym}) and{\samepage
\[
\partial_{u^k}\left(\gamma_{lj}\frac{\mathbf{b}_l}{\mathbf{b}_j}\right)
+\big(\gamma_{lk}\gamma_{kj}-g^{jk}h_j^0h_k^0\big)
\frac{\mathbf{b}_l}{\mathbf{b}_j}=0,\quad j,k,l\ \ \mathrm{distinct},
\]
which follows from (\ref{eq:gaj}) and (\ref{eq:hoj}).}

Note that (\ref{eq:sym}) can be written as
\begin{gather*}
\mathbf{b}_j^2(\partial_{u^j}\gamma_{jk}+2\gamma_{jk}\gamma_{kj})
+\mathbf{b}_k^2(\partial_{u^k}\gamma_{kj}+2\gamma_{kj}\gamma_{jk})
\nonumber\\
\qquad{}+\sum_{l\neq
j,k}\mathbf{b}_l^2\big(\partial_{u^k}\gamma_{lj}+2(\gamma_{lk}\gamma_{kj}
+\gamma_{lj}\gamma_{jk}-\gamma_{lj}\gamma_{lk})\big)=0,\quad j\neq
k,
\end{gather*}
so the dif\/ferential system (\ref{eq:hoj}), (\ref{eq:ajj}) is in
involution (completely integrable) for
\begin{gather}
\partial_{u^j}\gamma_{jk}=\partial_{u^k}\gamma_{kj}=-2\gamma_{jk}\gamma_{kj},\quad \ \
\partial_{u^k}\gamma_{lj}=2(\gamma_{lj}\gamma_{lk}-\gamma_{lk}\gamma_{kj}
-\gamma_{lj}\gamma_{jk}),\quad  j, k, l \ \
\mathrm{distinct}.\!\!\!\label{eq:comint}
\end{gather}
Thus if (\ref{eq:comint}) holds, then the solution of
(\ref{eq:hoj}), (\ref{eq:ajj}) is obtained by integrating $n$
f\/irst order ODE's (namely f\/inding the functions of $u^j$ upon
whose multiplication with $\mathbf{a}_j$ depends), so the space of
solutions depends on $(n-1)$ constants (the prime integral
$\sum_j\mathbf{b}_j^2+1=0$ removes a constant from the space of
solutions); if (\ref{eq:comint}) does not hold, then the space of
solutions depends on less than $(n-1)$ constants. Since for
Peterson's deformations of higher dimensional quadrics (or more
generally for deformations of sub-manifolds of the type
(\ref{eq:surf3}) with $f_k$, $g_k$, $h$ given by (\ref{eq:fgr}) with
$g_k(u^k)=0$, $k=2,\dots ,n-1$) we already have an $(n-1)$-dimensional
explicit family of deformations, we conclude that this family is
maximal and that (\ref{eq:comint}) holds in these cases.

\begin{remark}
Note that (\ref{eq:comint}) generalizes the case $n=2$ condition
$\partial_{u^1}\gamma_{12}=\partial_{u^2}\gamma_{21}=-2\gamma_{12}\gamma_{21}$
that the conjugate system $(u^1,u^2)$ is common to a Peterson's
$1$-dimensional family of deformations of surfaces (see Bianchi
\cite[Vol.~2, \S\S~294, 295]{B1}), so conjugate systems of
$n$-dimensional sub-manifolds in $\mathbb{C}^{n+1}$ satisfying
(\ref{eq:jkl}) and (\ref{eq:comint}) are a natural generalization
of Peterson's approach in the deformation problem.
\end{remark}

\section*{Acknowledgements}

I would like to thank the referees
for useful suggestions. The research has been supported by the University of Bucharest.

\pdfbookmark[1]{References}{ref}
\LastPageEnding


\begin{thebibliography}{99}

\footnotesize\itemsep=0pt

\bibitem{B}
Berger E., Bryant R.L., Grif\/f\/iths P.A.,
The Gauss equations and rigidity of isometric embeddings,
\href{http://dx.doi.org/10.1215/S0012-7094-83-05039-1}{{\it Duke Math. J.}} {\bf 50} (1983), 803--892.

\bibitem{B1}
Bianchi  L.,
Lezioni di geometria dif\/ferenziale, Vols.~1--4, Nicola Zanichelli Editore, Bologna, 1922, 1923, 1924, 1927.

\bibitem{Ca}
Calapso  P.,
Intorno alle superf\/icie applicabili sulle quadriche ed alle loro transformazioni,
\href{http://dx.doi.org/10.1007/BF02419392}{{\it Annali di Mat.}} {\bf 19}  (1912), no.~1, 61--82.\\
Calapso  P.,
Intorno alle superf\/icie applicabili sulle quadriche ed alle loro transformazioni,
\href{http://dx.doi.org/10.1007/BF02419394}{{\it Annali di Mat.}} {\bf 19}  (1912), no.~1, 107--157.

\bibitem{C}
Cartan \'E.,
Sur les vari\'{e}t\'{e}s de courboure constante d'un espace euclidien ou non-euclidien,
\href{http://www.numdam.org/item?id=BSMF_1919__47__125_1}{{\it Bull. Soc. Math. France}} {\bf 47} (1919), 125--160.\\
Cartan \'E.,
Sur les vari\'{e}t\'{e}s de courboure constante d'un espace euclidien ou non-euclidien,
\href{http://www.numdam.org/item?id=BSMF_1920__48__132_1}{{\it Bull. Soc. Math. France}} {\bf 48} (1920), 132--208.

\bibitem{D}
Darboux G.,
Le\c{c}ons sur la th\'{e}orie g\'{e}n\'{e}rale des surfaces et les applications g\'eom\'etriques du calcul inf\/init\'esimal, Vols.~1--4, Gauthier-Villars, Paris, 1894--1917.

\bibitem{H}
Horn  R.A., Johnson C.R.,
Matrix analysis, Cambridge University Press, Cambridge, 1985.

\bibitem{P} Peterson  K.-M., Sur la d\'{e}formation des surfaces du second ordre,
\href{http://afst.cedram.org/afst-bin/item?id=AFST_1905_2_7_1_69_0}{{\it Ann. Fac.  Sci. Toulouse S\'{e}r.~2\/}} {\bf 7}  (1905), no.~1, 69--107.

\end{thebibliography}
\end{document}